\documentclass[12pt]{article}%
\usepackage{amsfonts}
\usepackage{eurosym}
\usepackage{amssymb, amsthm}
\usepackage{amsmath}
\usepackage{amssymb}
\usepackage[allcolors=black,hidelinks]{hyperref}
\usepackage{graphicx}%
\providecommand{\U}[1]{\protect\rule{.1in}{.1in}}
\theoremstyle{plain}

\begin{document}

\begin{center}
{\LARGE \textbf{Team Problems and Stochastic\medskip}}

{\LARGE \textbf{ Programming}}

\bigskip\bigskip

Igor V. Evstigneev\footnote{Department of Economics, University of Manchester,
Oxford Road, Manchester, M13 9PL, UK.\ E-mail::
igor.evstigneev@manchester.ac.uk.}, Mohammad J. Vanaei\footnote{Economics
Department, University of Manchester, Oxford Road, Manchester, M13 9PL UK.
E-mail: mohammadjavad.vanaei@manchester.ac.uk\medskip.}, and Mikhail V.
Zhitlukhin\footnote{Steklov Mathematical Institute, Russian Academy of
Sciences, Moscow, 119991, Russia. E-mail: mikhailzh@mi-ras.ru. (Corresponding
author.)}\bigskip\bigskip

\bigskip
\end{center}

\begin{quotation}
\textbf{Abstract.} The paper bridges two vast areas of research: stochastic
team decision problems and convex stochastic programming. New methods
developed in the latter are applied to the study of fundamental problems in
the former. The main results are concerned with the Lagrangian relaxation of
informational and material constraints in convex stochastic team
problems.\bigskip

\textbf{Keywords} \ Team decision problems, Stochastic programming,
Informational constraints, Material constraints, Lagrangian relaxation

\textbf{Mathematics Subject Classification -- MSC2020:} 93E20, 90C15, 90C46,
90B70 \bigskip\ 
\end{quotation}

\section{Introduction}

\textbf{1.} \textbf{Team problems.} In a stochastic team problem, multiple
decision-makers (DMs) or agents operate within a system influenced by random
factors. Each DM observes some information about the system and chooses an
action. The combined actions of all DMs affect the overall system state and a
global objective function. The goal lies in coordinating these actions to
achieve the best possible outcome, often under informational or material
constraints. Informational constraints involve the limited knowledge or
incomplete information available to the decision-maker. Material constraints
typically refer to resource limitations that restrict choices, such as budget,
time, or physical capabilities. These constraints shape how individuals make
decisions, even when they are assumed to be rational.

The theory of decision making in teams (with the emphasis on economic
applications) was founded in the works of Marschak \cite{Marschak1955},
Marschak and Radner \cite{MarschakRadner1972}, and Radner
\cite{Radner1959,Radner1962,Radner1972,Radner1972a}. For early contributions
to the field see also Witsenhausen \cite{Witsenhausen1968}, Ho and Chu
\cite{HoChu1972}, Groves \cite{Groves1973}, Sandell et al.
\cite{Sandell-et-al1978}, Arrow and Radner \cite{ArrowRadner1979}, and Ho
\cite{Ho1980}.\ Due to the conceptional and technical difficulties, the theory
has been slow to develop; for a discussion of key problems in the area see
Papadimitriou and Tsitsiklis \cite{PapadimitriouTsitsiklis1986}. Substantial
progress has been made during the last two decades. Comprehensive reviews of
the current state of the art of the field can be found in the recent papers by
Zhang et al. \cite{Zhang-et-al-2020}, Malikopoulos \cite{Malikopoulos2023} and
in the monographs by Y\"{u}ksel and Ba\c{s}ar
\cite{YukselBasar2013,YukselBasar2024}.

\textbf{2. Stochastic programming.} In stochastic programming problems, a
single DM makes decisions at consecutive moments of time observing events from
an increasing sequence of $\sigma$-algebras in the underlying probability
space and using this information in decision making. The goal consists in the
maximization of a given functional, typically having the form of the
expectation of an objective function depending on the whole sequence of
decisions. A characteristic feature of stochastic programming problems (in
contrast, say, to Markov decision processes) is that the decisions made do not
influence the probability distribution on the given probability space.
Although the realizations of random events are not known in advance, the DM
knows their probabilities, as well as the probability distributions of
observable random variables.

An important role in the development of the mathematical theory of stochastic
programming was played by the seminal work of Dynkin \cite{Dynkin1972}, Dynkin
and Yushkevich \cite[Chapter 9]{DynkinYushkevich1979}\footnote{The main
results of \cite{Dynkin1972} are presented in Chapter 9 of the book
\cite{DynkinYushkevich1979}.}, Rockafellar \cite{Rockafellar1974}, and
Rockafellar and Wets
\cite{RockafellarWets1976b,RockafellarWets1976,RockafellarWets1976a,
RockafellarWets1976c,RockafellarWets1976c,RockafellarWets1976c,RockafellarWets1978}%
. These studies served as the starting point for that line of research on
stochastic programming to which this paper pertains. Contributions to the
field were made, among others, by Evstigneev
\cite{Evstigneev1976,Evstigneev1976a} and Pennanen and Perkki\"{o}
\cite{PennanenPerkkio2012,PennanenPerkkio2018,PennanenPerkkio2023,PennanenPerkkio2024,PennanenPerkkio2025}%
. Comprehensive reviews of achievements in this area are provided in Arkin and
Evstigneev \cite{ArkinEvstigneev1987}, Shapiro et al.
\cite{ShapiroDentchevaRuszczynski2021} and Pennanen and Perkki\"{o}
\cite{PennanenPerkkio2024}, where the reader can find further references. Most
of the results in the above literature are motivated by the economic and
financial applications.

\textbf{3. Integration of team problems and stochastic programming.} In this
paper we study team decision problems in a framework extending conventional
stochastic programming. The main distinction of the present setting from the
conventional ones is that at each moment of time there are several (rather
than one) DMs, making decisions simultaneously and independently. Furthermore,
there is a partial ordering on the set of all DMs describing interrelations
and interdependencies between the agents in space and time. The informational
structure is represented by a family of $\sigma$-algebras in the underlying
probability space, describing information available for each agent in the
team. The ordering (by inclusion) between these $\sigma$-algebras is assumed
to be compatible with the ordering given on the set of DMs. Thus, the
optimization problems we study here may be characterized in a nutshell as
stochastic programming problems with "partially ordered time", cf. Mandelbaum
and Vanderbei \cite{MandelbaumVanderbei1981}, Lawler and Vanderbei
\cite{LawlerVanderbei1983}, and Evstigneev and Zhitlukhin
\cite{EvstigneevZhitlukhin2013}.

This paper focuses on Lagrangian relaxation of informational and material
constraints in team problems. A key role in this work is played by the
Rockafellar-Wets approach to the analysis of information constraints in convex
stochastic optimization problems, see Rockafellar \cite{Rockafellar1974} and
Rockafellar and Wets \cite{RockafellarWets1976a}. The idea of the approach
lies in the representation of measurability conditions, which is a standard
way of describing informational restrictions, as linear operator constraints
in spaces of measurable functions. With a finite probability space and under
the assumptions of convexity, this makes the problem computationally tractable
thanks to the existence of efficient algorithms for numerical solutions of
convex problems with linear constraints. Following the Rockafellar-Wets
approach, we construct Lagrange multipliers associated with information
constraints ("shadow prices on information"), which not only represent a
useful theoretical tool, but also have important applications in economics and
finance---see, e.g., Evstigneev et al.
\cite{EvstigneevKleinHaneveldandMirman1999} and Fl\aa m and Koutsougeras
\cite{FlamKoutsougeras2010}. For related results on information and duality in
convex stochastic optimization see Pennanen and Perkki\"{o}
\cite{PennanenPerkkio2018}.

The paper is organized as follows. Section 2 describes the model and states
the results. Section 3 gives their proofs. Section 4 states corollaries to the
main results and discusses applications. The Appendix assembles several facts
from functional analysis used in this work.

\section{Model and results}

\textbf{1.} \textbf{The model.} There is a finite set $T$ of decision makers
(DMs) or agents acting at moments of time $n=0,1,...,N$. We denote by $n(t)$
the moment of time when DM $t\in T$ operates and by $T_{n}:=\{t:n(t)=n\}$ the
set of agents making decisions at time $n$. It is assumed that $T_{n}%
\neq\emptyset$ for each $n=0,1,...,N$. The set $T$ is endowed with a partial
ordering $\preccurlyeq$ describing interrelations and interdependences (in
space and time) between DMs $t\in T$. We write $t\prec s$ if $t\preccurlyeq s$
and $t\neq s$.

It is supposed that if $t\prec s$, then $n(t)<n(s)$. If $t,s\in T_{n}$ for
some $n$, then $t$ and $s$ are incomparable in terms of the strict partial
ordering $\prec$: neither of the relations $t\prec s$ and $s\prec t$ holds.
This means that agents acting at the same moment of time make decisions
simultaneously and independently. The set $T\ $with the strict partial order
$\prec$ can equivalently be thought of as an acyclic directed graph such that
$t\prec s$ if and only if there is a path from "past" to "future" along the
edges of the graph leading from $t$ to $s$.

To each $t$ there corresponds a $\sigma$-algebra $\mathcal{F}_{t}%
\subseteq\mathcal{F}$ in the probability space $(\Omega,\mathcal{F},P)$
interpreted as the class of events observable by agent $t$ and affecting $t$'s
decisions. It is supposed that if $t\prec s$, then $\mathcal{F}_{t}%
\subseteq\mathcal{F}_{s}$. Decision variables of DMs are represented by
$\mathcal{F}_{t}$-measurable functions $x_{t}:\Omega\rightarrow V_{t}$, where
$V_{t}$ is an $m_{t}$-dimensional vector space . The goal of the team of DMs
lies in the maximization of a common functional $F(x)$ of $x=(x_{t}%
)_{t\in\mathbb{T}}$ subject to informational and material constraints.

For each $t\in T$ and $\omega$, let $A_{t}(\omega)$ be a subset in the space
$V_{t}$. Elements $a\in A_{t}(\omega)$ represent admissible decisions of agent
$t$ in the random situation $\omega$. A collection $x=(x_{t})_{t\in T}$ of
$\mathcal{F}$-measurable vector functions $x_{t}:\Omega\rightarrow V_{t}$ will
be called a \textit{program} if%
\begin{equation}
x_{t}(\omega)\in A_{t}(\omega)\text{ a.s.},\text{ }t\in T. \label{f1}%
\end{equation}
("A.s." means "almost surely" with respect to the probability $P$.) The class
of all programs will be denoted by $X$. They describe all possible decisions
of the team using full information about $\omega$.

Suppose that for each $t$ we are given a function $g_{t}(\omega,a^{t})$, where
$a^{t}=(a_{i})_{i\preccurlyeq t}$,$\ a_{i}\in V_{i}$, with values in an
$l_{t}$-dimensional vector space $W_{t}$. Let $F(x)$ be a real-valued
functional defined on the set $X$ of all programs. The following stochastic
optimization problem will be studied.

($\mathcal{P}$) Maximize $F(x)$ among all programs $x=(x_{t})_{t\in T}\in X$
satisfying for all $t$ the following conditions:%

\begin{equation}
g_{t}(\omega,x^{t}(\omega))\geq0\ \text{a.s.},\ \label{f2}%
\end{equation}

\begin{equation}
x_{t}(\omega)\text{ is }\mathcal{F}_{t}\text{-measurable.} \label{f3}%
\end{equation}

Thus, each DM $t$ for each random situation $\omega$ has to choose a decision
$x_{t}(\omega)$ so that constraints (\ref{f2}) and (\ref{f3}) are met, and the
program $(x_{t})$ maximizes the functional $F$. Inequalities (\ref{f2})
represent \textit{material constraints} expressing in economic problems, e.g.,
resource, capacity, manpower, or other limitations. Requirement (\ref{f3})
means that the choice of $x_{t}(\omega)$ depends only on events observable by
agent $t$ , so that (\ref{f3}) may be regarded as an \textit{information
constraint}. Programs satisfying the information constraints (\ref{f3}) will
be called \textit{adapted} (to the system of $\sigma$-algebras $\mathcal{F}%
_{t}$). We will denote their totality by $H$. Adapted programs for which
condition (\ref{f2}) holds will be called \textit{feasible}. Problem
($\mathcal{P}$) consists in finding an \textit{optimal }feasible program---the
one maximizing the functional $F(x)$ on the set of all feasible programs.

\textbf{2. Stochastic Lagrange multipliers.} Assume that the following
conditions (A) and (B) are fulfilled.

\begin{itemize}
\item[(A)] For each $t$, the set $A_{t}(\omega)$ is convex, closed, uniformly
bounded and depends $\mathcal{F}_{t}$-measurably on $\omega$. The last
condition means that for each point $a$ the distance between $a$ and
$A_{t}(\omega)$ is an $\mathcal{F}_{t}$-measurable function of $\omega$.
\end{itemize}

Put%
\[
V^{t}=%
{\textstyle\prod_{i\preccurlyeq t}}
V_{i},\ \mathcal{B}^{t}=\mathcal{B}(V^{t})=%
{\textstyle\prod_{i\preccurlyeq t}}
\mathcal{B}(V_{i}),,
\]
where $\mathcal{B}(\cdot)$ is the Borel $\sigma$-algebra.

\begin{itemize}
\item[(B)] The functions $g_{t}(\omega,a^{t})$ are measurable with respect to
$\mathcal{F}_{t}\times\mathcal{B}^{t}$. For every $t$ and $\omega$, the
function $g_{t}(\omega,a^{t})$ is continuous and concave (coordinate-wise)
with respect to $a^{t}\in A^{t}(\omega):=%
{\textstyle\prod_{i\preccurlyeq t}}
A_{i}(\omega)$. There exists a number $c$ such that $|g_{t}(\omega,a^{t})|\leq
c$ for all $\omega,t,a^{t}$.
\end{itemize}

Our main results are concerned with the construction of stochastic Lagrange
multipliers relaxing material and informational constraints in problem
($\mathcal{P}$). It will be convenient to represent the information constraint
(\ref{f3}) in an equivalent form:
\begin{equation}
x_{t}-E_{t}x_{t}=0, \label{f3-prime}%
\end{equation}
where $E_{t}x_{t}=E(x_{t}\mid\mathcal{F}_{t})$ is the conditional expectation
given the $\sigma$-algebra $\mathcal{F}_{t}$ with respect to the probability
$P$.

\textbf{Definition.} We say that random functions $p_{t}(\omega)$ and
$q_{t}(\omega)$ are \textit{stochastic Lagrange multipliers relaxing
constraints (\ref{f2}) and (\ref{f3-prime}) }in the optimization problem
($\mathcal{P}$) if the following conditions hold:

\begin{itemize}
\item[(a)] $p_{t}(\omega)$ is a non-negative $\mathcal{F}_{t}$-measurable
integrable function with values in $W_{t}$.

\item[(b)] $q_{t}(\omega)$ is an $\mathcal{F}$-measurable integrable function
with values in $V_{t}$.

\item[(c)] For all programs $x=(x_{t})_{t\in T}\in X$ and each solution
$\bar{x}$ to problem ($\mathcal{P}$), the \textit{Kuhn-Tucker conditions}
hold:%
\begin{equation}
L(x)\leq L(\bar{x}), \label{L}%
\end{equation}
and%
\begin{equation}
Ep_{t}g_{t}(\bar{x}^{t})=0\text{, for each }t\in T, \label{CS}%
\end{equation}

\end{itemize}

where
\begin{equation}
L(x):=F(x)+\sum_{t\in T}Ep_{t}g_{t}(x^{t})+\sum_{t\in T}Eq_{t}(x_{t}%
-E_{t}x_{t}) \label{L-x}%
\end{equation}
is the \textit{Lagrangian} for the problem ($\mathcal{P}$). Here "$E$" stands
for the expectation with respect to the probability $P$. Throughout the paper,
we often omit the argument $\omega$ for the sake of brevity and write, e.g.,
$Ep_{t}g_{t}(x^{t})$ instead of $Ep_{t}(\omega)g_{t}(\omega,x^{t}(\omega))$.

Observe that for a feasible program $\bar{x}$ conditions (\ref{L})\&(\ref{CS})
hold if and only if%
\begin{equation}
F(x)+\sum_{t\in T}Ep_{t}g_{t}(x^{t})+\sum_{t\in T}Eq_{t}(x_{t}-E_{t}x_{t})\leq
F(\bar{x}) \label{f4}%
\end{equation}
for all programs $x=(x_{t})_{t\in T}\in X$. Indeed, if (\ref{f4}) is valid,
then by substituting the feasible program $\bar{x}$ into the left-hand side of
(\ref{f4}) and using the fact that it satisfies (\ref{f3-prime}), we get
$F(\bar{x})=L(\bar{x})$, which yields (\ref{L}) and (\ref{CS}). Conversely, if
(\ref{L})\&(\ref{CS}) hold and $\bar{x}$ is feasible, then $F(\bar{x}%
)=L(\bar{x})$ by virtue of (\ref{f3-prime}) and (\ref{CS}), and so (\ref{L})
implies (\ref{f4}).

We can see from the above considerations that if the Lagrange multipliers
$p_{t}$ and $q_{t}$ are known, then all solutions to ($\mathcal{P}$) can be
found as follows. It is sufficient to find all programs $x$ maximizing $L(x)$
on $X$ and choose among them those ones $\bar{x}$ that are feasible and
satisfy the \textit{complementary slackness condition }(\ref{CS}). Then
(\ref{f4}) will be valid as a consequence of (\ref{L})\&(\ref{CS}), which will
imply $F(x)\leq F(\bar{x})$ in view of (\ref{f2}) and (\ref{f3-prime}). Thus
(\ref{L})\&(\ref{CS}) are necessary and sufficient conditions for a feasible
program $\bar{x}$ to be optimal.

The idea of representing the measurability requirement (\ref{f3}) in the form
(\ref{f3-prime}) belongs to Rockafellar \cite{Rockafellar1974}. It turned out
to be quite fruitful. The representation of condition (\ref{f3}) as a linear
operator constraint in a space of measurable functions opened possibilities
for the application of the methods of Lagrangian relaxation to the analysis of
information constraints. Initially, this idea was developed by Rockafellar and
Wets \cite{RockafellarWets1976a} and then by many others, see the monograph by
Pennanen and Perkki\"{o} \cite{PennanenPerkkio2024} and references therein.
For a discussion of applications see Section 4.

\textbf{3. The main result.} We now formulate the assumptions under which the
main result (Theorem 1) holds. Recall that $X$ denotes the set of all programs
and $H$ the class of all adapted programs.

\begin{itemize}
\item[(C)] The functional $F$ is concave on $X$ and continuous on $X$ with
respect to a.s. convergence.

\item[(D)] There exists an adapted program $(\mathring{x}_{t})\in H$ and a
positive real number $\theta$ such that, for any $t$ and $\omega$, the set
$A_{t}(\omega)$ contains a ball of the radius $\theta$ with the centre at
$\mathring{x}_{t}(\omega)$. (Slater's condition for constraints (\ref{f3})).

\item[(E)] There exists an adapted program$\ \hat{x}=(\hat{x}_{t})\in H$ and a
constant vector $\kappa>0$ for which $g_{t}(\hat{x}^{t})\geq\kappa$ ($t\in
T$). (Slater's condition for constraints (\ref{f2})).
\end{itemize}

All inequalities between vectors---strict and non-strict---are understood
coordinate-wise. Let $x=(x_{t})$ be an adapted program and let $s$ be and
element in $T$. We will say that $x$ is $s$\textit{-feasible} if it satisfies
constraints (\ref{f2}) for all $t\prec s$.

\begin{itemize}
\item[(F)] For every $s\in T$ and each $s$-feasible adapted program
$(x_{t})\in H$, there exists $(\tilde{x}_{t})\in H$ such that $x_{t}=\tilde
{x}_{t}$ for $t\neq s$, and $g_{s}(\tilde{x}^{s})\geq0$ (a.s.).
\end{itemize}

This condition means that it is possible to alter ("switch to a safe mode") an
$s$-feasible program $(x_{t})$ at the node $s$ of the graph $T$ so that
constraint (\ref{f2}) is satisfied for $t=s$. Assumption (F) is akin to the
condition of \textit{relatively complete recourse}, well known in stochastic
programming, see, e.g., Rockafellar and Wets \cite{RockafellarWets1976}.

\textbf{Theorem 1. }\textit{If conditions (A)-(F) are satisfied, then the
maximum in problem (}$\mathcal{P}$\textit{) is attained and there exist
stochastic Lagrange multipliers }$p_{t}$\textit{ and }$q_{t}$\textit{ relaxing
constraints (\ref{f2}) and (\ref{f3-prime}).}

Assumptions (A)-(F) will be supposed to hold in what follows.\textit{ }

\textbf{4. Constraint relaxation almost surely.} We shall now consider the
special case when $F$ is an integral functional. Let $f(\omega,a)$ be a
real-valued function defined for $\omega\in\Omega$ and
\[
a=(a_{t})_{t\in T}\in A(\omega):=\prod_{t\in T}A_{t}(\omega)
\]
Suppose that the following conditions hold:

\begin{itemize}
\item[(C1)] The function
\[
f^{\infty}(\omega,a)=\left\{
\begin{array}
[c]{cc}%
f(\omega,a), & a\in A(\omega),\\
\infty, & a\notin A(\omega),
\end{array}
\right.
\]
is measurable with respect to $(\omega,a)$.

\item[(C2)] The function $f(\omega,a)$ is concave and upper semi-continuous in
$a\in A(\omega)$.

\item[(C3)] There is an integrable function $d(\omega)$ such that
$|f(\omega,a)|\leq d(\omega)$ for $a\in A(\omega)$ and $\omega\in\Omega$.
\end{itemize}

Let us define the functional $F$ of $x=(x_{t})_{t\in T}\in X$ by the formula
\begin{equation}
F(x)=Ef(\omega,x(\omega)). \label{E-f}%
\end{equation}
Then condition (C) is satisfied, and the Lagrangian for the problem
($\mathcal{P}$) can be represented in the form
\begin{equation}
L(x)=El(\omega,x(\omega)), \label{E-l}%
\end{equation}
where%
\begin{equation}
l(\omega,a):=f(\omega,a)+\sum_{t=1}^{N}\left[  p_{t}(\omega)g_{t}(\omega
,a^{t})+\left(  q_{t}(\omega)-\tilde{q}_{t}(\omega)\right)  a_{t}\right]  ,
\label{l}%
\end{equation}
and
\begin{equation}
\tilde{q}_{t}(\omega):=E_{t}q_{t}(\omega). \label{q-tilde}%
\end{equation}
We use here the fact that%
\[
Eq_{t}(x_{t}-E_{t}x_{t})=Eq_{t}x_{t}-Eq_{t}E_{t}x_{t}=Eq_{t}x_{t}-E(E_{t}%
q_{t})x_{t}=E(q_{t}-\tilde{q}_{t})x_{t}.
\]

\textbf{Proposition 1.} \textit{Suppose that the functional }$F(x)$\textit{ is
of the form (\ref{E-f}), where the function }$f$\textit{\ satisfies conditions
(C1)-(C3). Then functions }$p_{t}$\textit{ and }$q_{t}$\textit{ satisfying (a)
and (b) represent stochastic Lagrange multipliers }$p_{t}$\textit{ and }%
$q_{t}$\textit{ relaxing constraints (\ref{f2}) and (\ref{f3-prime}) if and
only if }%
\begin{equation}
\max_{a\in A(\omega)}l(\omega,a)=l(\omega,\bar{x}(\omega))\mathit{\ }%
\text{\textit{a.s.}}, \label{l-a-s}%
\end{equation}%
\begin{equation}
p_{t}g_{t}(\bar{x}^{t})=0\text{\textit{ a.s. for all }}t\in T. \label{CS-a-s}%
\end{equation}

Relations (\ref{l-a-s}) (\ref{CS-a-s}) represent almost sure versions of the
Kuhn-Tucker conditions (\ref{L}) and (\ref{CS}). If they hold, we say that the
Lagrange multipliers $p_{t}$ and $q_{t}$ \textit{relax constraints (\ref{f2})
and (\ref{f3-prime}) almost surely}. If $p_{t}$ and $q_{t}$ are known, then
the original (generally, infinite-dimensional) convex optimization problem
essentially reduces to a family of finite-dimensional ones, which are solvable
point-wise and amenable to efficient computational algorithms.

By combining Theorem 1 and Proposition 1, we obtain

\textbf{Theorem 2.} \textit{Under assumptions (C1)-(C3), there exist
stochastic Lagrange multipliers }$p_{t}$\textit{ and }$q_{t}$\textit{ relaxing
constraints (\ref{f2}) and (\ref{f3-prime}) almost surely.}

Recall that we always assume that conditions (A)-(F) hold.\textit{ }

\textbf{Remark.} Note that there is a broad range of applications, where the
functional $F(x)$ satisfies (C) but is of a more general type than
(\ref{E-l}). Non-integral functionals appear, in particular, in connection
with risk measures, see e.g. Szeg\"{o} \cite{Szego2004}. An important example
of this kind is the Markowitz \cite{Markowitz1959} mean-variance model of
portfolio selection, where the role of a risk measure is played by the
variance of a random variable.

\section{Proofs}

\textbf{1. The maximum in (}$\mathcal{P}$\textbf{) is attained. }We begin with
the \textit{proof of Theorem 1}. The first step is to show that a solution to
problem ($\mathcal{P}$) exists. Denote by $Q$ the set of all vector functions%
\[
x=(x_{t})_{t\in T}\in L_{1}(\mathcal{F},V),\text{ }V:=\prod_{t\in T}V_{t},
\]
such that $x=x^{\prime}$ a.s. for some $x^{\prime}\in X$ satisfying
constraints (\ref{f2}) and (\ref{f3}). Clearly, $Q$ is convex and closed with
respect to a.s. convergence. Hence, by (FA4)\footnote{For asserions (F1) --
(F8) see the Appendix.}, there exists an $x\in Q$ which maximizes $F$ on $Q$.
By the definition of $Q$, there is a sequence $\bar{x}$ which coincides with
$x$ a.s. and satisfies the constraints (\ref{f1})--(\ref{f3}). This sequence
$\bar{x}$ is a solution to problem ($\mathcal{P}$).

The proof of the existence of $p_{t}$ and $q_{t}$ will be divided into 4 steps.

\textbf{2. Relaxation of material constraints in }$(L_{\infty})^{\ast}$. We
first construct non-negative functionals $\pi_{t}\in\left[  L_{\infty}\left(
\mathcal{F}_{t},W_{t}\right)  \right]  ^{\ast}$ with the property
\begin{equation}
F(x)+\sum_{t\in T}\pi_{t}g_{t}(x^{t})\leq F(\bar{x})\quad(x\in H). \label{f9}%
\end{equation}
Let us regard an adapted program $x=(x_{t})\in H$ as an element of the space
$D_{1}=\prod_{t\in T}L_{\infty}\left(  \mathcal{F}_{t},V_{t}\right)  $, and
$g(x)=(g_{t}(x^{t}))_{t\in T}$ as a mapping of $D_{1}$ into the space
$D_{2}=\prod_{t\in T}L_{\infty}\left(  \mathcal{F}_{t},W_{t}\right)  $ ordered
in the canonical way. In view of the assumption (D), the proposition (FA2) can
be applied to the functional $F$, the set $H$ and the mapping $g:H\rightarrow
D_{2}$. According to this proposition, there exists a non-negative functional
$\pi=(\pi_{t})_{t\in T}\in D_{2}^{\ast}$ satisfying the condition (\ref{f9}).

\textbf{3. Relaxation of material constraints in }$L_{1}$.We shall now prove
that if the functionals $\pi_{t}$ have the property (\ref{f9}), then their
absolutely continuous components $\pi_{t}^{a}$ (see (FA5)) have the same property.

Denote by $I$ the number of elements in the set $T$. Let us enumerate all
elements $t_{1},...,t_{I}$ in the set $T$ in such a way that if $t_{i}\prec
t_{m}$, then $i<m$. This can be done as follows. First enumerate all elements
in $T_{1}$ in any order, then in $T_{2}$, then in $T_{3}$, and so on. We
assumed (see \textbf{2.1}) that $t_{i}\prec t_{m}$ implies $n(t_{i})<n(t_{m}%
)$. Consequently, $i<m$ by the construction of the sequence $t_{1},...,t_{I}$.

Let us write for shortness $i$ in place of $t_{i}$. Nodes of the graph $T$
will be identified with natural numbers $i=1,2,...,I$ and the relation
$t_{i}\prec t_{m}$ will be written as $i\prec m$. As we have shown above,
\begin{equation}
i\prec m\ \Longrightarrow i<m. \label{order}%
\end{equation}
A program $x=(x_{t})_{t\in T}$ can be equivalently represented as%
\[
x=(x_{t})_{t\in T}=(x_{t_{i}})_{i=1}^{I}=(x_{i})_{i=1}^{I}.
\]
Fix a natural number $m=1,...,I$ and an adapted program $x=(x_{i})_{i=1}%
^{I}\in H$. Consider the inequality
\begin{equation}
F(x)+\pi^{a}g(x)+\sum_{i=1}^{m}\pi_{i}^{s}g_{i}(x^{i})\leq F(\bar{x}%
)\quad(x\in H), \label{f10}%
\end{equation}
where where $\pi^{a}=(\pi_{i}^{a})_{i=1}^{I}$. By virtue of (\ref{f9}), this
inequality is valid for $m=I$ \ and all adapted programs $x=(x_{i})_{i=1}%
^{I}\in H$. We have to verify it for $m=0$ (when $\sum_{i=1}^{0}=0$ by
definition). To this end it is sufficient to prove the following fact. If
inequality (\ref{f10}) is true for some $m\geq1$, then it is true for $m-1$.

Fix an $m\geq1$ and an adapted program $x=(x_{i})_{i=1}^{I}\in H$. By virtue
of condition (F), there exists an adapted program $(\tilde{x}_{i})_{i=1}%
^{M}\in H$ such that $\tilde{x}_{i}=x_{i}$ for $i\neq m$, and $g_{m}(\tilde
{x}^{m})\geq0$ (a.s.). Let $\Gamma_{k}$, $k=1,2,...$, be the sets
corresponding to the singular functional $\pi_{m}^{s}$ (see (FA5). For each
$k$ consider the adapted program $\tilde{x}=(\tilde{x}_{i}^{k})_{i=1}^{I}$
defined by
\[
\tilde{x}_{i}^{k}=\left\{
\begin{array}
[c]{cc}%
x_{i} & \text{if }i\neq m,\\
\chi_{\Gamma_{k}}\tilde{x}_{m}+\chi_{\bar{\Gamma}_{k}}x_{m} & \text{if\ }i=m.
\end{array}
\right.
\]
By substituting $(\tilde{x}_{i}^{k})_{i=1}^{I}$ into (\ref{f10}), we get
\[
F(x^{k})+\pi^{a}g(x^{k})+\sum_{i=1}^{m-1}\pi_{i}^{s}g_{i}((\tilde{x}^{k}%
)^{i})+\pi_{m}^{s}\left(  g_{m}(x^{m})\chi_{\bar{\Gamma}_{k}}\right)
\]%
\[
+\pi_{m}^{s}\left(  g_{m}(\tilde{x}^{m})\chi_{\Gamma_{k}}\right)  \leq
F(\bar{x}).
\]
The last term of the left-hand side of this inequality is non-negative by the
definition of the program $\tilde{x}$ and the term $\pi_{m}^{s}\left(
g_{m}(x^{m})\chi_{\bar{\Gamma}_{k}}\right)  $ is equal to $0$ by the
definition of the sets $\Gamma_{k}$. Hence these terms can be omitted, which
yields%
\begin{equation}
F(x^{k})+\pi^{a}g(x^{k})+\sum_{i=1}^{m-1}\pi_{i}^{s}g_{i}((\tilde{x}^{k}%
)^{i})\leq F(\bar{x}). \label{sum-m-1}%
\end{equation}
Observe that if $i\leq m-1$ (as in the left-hand side of (\ref{sum-m-1})), or
equivalently, $i<m$, then
\begin{equation}
(\tilde{x}^{k})^{i}=(\tilde{x}_{j}^{k})_{j\preccurlyeq i}=(x_{j}%
)_{j\preccurlyeq i}=x^{i} \label{x-k-i}%
\end{equation}
because (a) $\tilde{x}_{j}^{k}=x_{j}$ for $j\neq m$, and (b)\ the set
$\{j:j\preccurlyeq i\}$ does not contain $m$. Indeed, suppose the contrary:
$m\preccurlyeq i$. Since $i<m$, we have $i\neq m$, and so $m\prec i$. Then by
virtue of (\ref{order}), $m<i$, which contradicts our assumption.

By combining (\ref{sum-m-1}) and (\ref{x-k-i}), we write%
\[
F(x^{k})+\pi^{a}g(x^{k})+\sum_{i=1}^{m-1}\pi_{i}^{s}g_{i}(x^{i})\leq F(\bar
{x}).
\]
Now it remains to pass to the limit as $k\rightarrow\infty$. In view of (C),
this will yield (\ref{f10}) with $m$ replaced by $m-1$.

\textbf{4. Relaxation of informational constraints in }$(L_{\infty})^{\ast}%
$\textbf{.} Let $p_{t}\in L_{1}(\mathcal{F}_{t},W_{t})$ be the functions
corresponding to the functionals $\pi_{t}^{a}$. It has been just established
that
\begin{equation}
F(x)+\sum Ep_{t}g_{t}(x^{t})\leq F(\bar{x})\quad(x\in H). \label{f11}%
\end{equation}
We shall now construct $q_{t}\in L_{\infty}(\mathcal{F},V_{t})$ ($t=1,2,\ldots
,N$) such that
\begin{equation}
F^{\prime}(x)+\sum Eq_{t}\left(  x_{t}-E(x_{t}\mid\mathcal{F}_{t})\right)
\leq F^{\prime}(\bar{x})\quad(x\in X), \label{f12}%
\end{equation}
where $F^{\prime}(x)$ is the expression at the left-hand side of (\ref{f11}).
Theorem 1 will then follow immediately, since (\ref{f11}) and (\ref{f12})
imply (C).

For any random variable $\xi$, let us write $E_{t}\xi$ instead of
$E(\xi|\mathcal{F}_{t})$ and denote $\xi-E_{t}\xi$ by $h_{t}(\xi)$. Our aim
now is to find a continuous linear functional $\rho=(\rho_{1},\ldots,\rho
_{N})$ on the space $D=L_{\infty}(\mathcal{F},V_{1})\times\cdots\times
L_{\infty}(\mathcal{F},V_{N})$ such that
\begin{equation}
F^{\prime}(x)+\sum\rho_{t}h_{t}(x_{t})\leq F^{\prime}(\bar{x})\quad(x\in X).
\label{f13}%
\end{equation}
Observe that the constraints (\ref{f3}) can be rewritten in the form $h(x)=0$,
where $h(x)=\{h_{1}(x_{1}),\ldots,h_{N}(x_{N})\}$. The linear operator
$h:D\rightarrow D$ is continuous and maps $D$ onto its closed subspace
$\{(x_{1},\ldots,x_{N}):E_{t}x_{t}=0,\,t=1,2,\ldots,N\}$. Furthermore,
$\{x:h(x)=0\}\cap\operatorname{Int}X\neq\emptyset$ by virtue of the condition
(E). Thus, according to (FA1), there exists a functional $\rho$ with the
property (\ref{f13}).

\textbf{5. Relaxation of informational constraints in }$L_{1}$\textbf{.} Let
us verify that (\ref{f13}) remains valid if we replace $\rho_{t}$ by $\rho
_{t}^{a}$. This will complete the proof of Theorem 1. Fix $t=m$, $x\in X$ and
a sequence of sets $\Gamma_{k}$ which is \textquotedblleft
charged\textquotedblright\ by the functional $\rho_{m}^{s}$ and has the
additional property (\ref{f8}). Using $\Gamma_{k}$, we define for each $k$ the
sequence $x^{[k]}=\{x_{t}^{[k]}\}$:
\begin{equation}
x_{t}^{[k]}=x_{t}\quad(t\neq m);\quad x_{m}^{[k]}=\chi_{\Gamma_{k}}E_{m}%
x_{m}+\chi_{\bar{\Gamma}_{k}}x_{m}. \label{f14}%
\end{equation}
This sequence belongs to the class $X$, since $E_{m}{x_{m}}\in A_{m}(\omega)$
a.s. (see (FA8)). Thus, by (\ref{f13}),
\begin{equation}
F^{\prime\lbrack k]})+\sum_{t\neq m}\rho_{t}(h_{t}({x_{t}}))+\rho_{m}%
^{a}(h_{m}(x_{m}^{[k]}))+\rho_{m}^{s}(h_{m}(x_{m}^{[k]}))\leq F^{\prime}%
(\bar{x}), \label{f14a}%
\end{equation}
and we need only to show that
\begin{equation}
r_{k}\equiv\rho_{m}^{s}(h_{m}(x_{m}^{[k]}))\rightarrow0\quad\text{as
}k\rightarrow\infty. \label{f14b}%
\end{equation}
This, in turn, will be proved if we establish that $\Vert\psi_{k}\Vert
_{\infty}\rightarrow0$, where $\psi_{k}=\chi_{\Gamma_{k}}h_{m}({x_{m}^{[k]}}%
)$. Indeed, $\rho_{m}^{s}$ is continuous with respect to $\Vert\cdot
\Vert_{\infty}$ and $r_{k}=\rho_{m}^{s}(\psi_{k})$ by the definition of
$\Gamma_{k}$. It follows from (\ref{f14}) that, for $\omega\in\Gamma_{k}$,
\begin{equation}
\psi_{k}=x_{m}^{[k]}-E_{m}x_{m}^{[k]}=E_{m}x_{m}-E_{m}(\chi_{\bar{\Gamma}_{k}%
}x_{m})-E_{m}(x_{m}E_{m}\chi_{\Gamma_{k}}). \label{f14c}%
\end{equation}
The last expression can be rewritten in the form
\begin{equation}
E_{m}\left(  \chi_{\bar{\Gamma}_{k}}E_{m}(\chi_{\Gamma_{k}}x_{m})\right)
-E_{m}\left(  x_{m}\chi_{\bar{\Gamma}_{k}}E_{m}\chi_{\Gamma_{k}}\right)  ,
\label{f14d}%
\end{equation}
and we have $\Vert\psi_{k}\Vert_{\infty}\rightarrow0$ by virtue of (\ref{f8}).
The proof of Theorem\ 1 is complete.\hfill$\square$

\textbf{6. Proof of Proposition 1}\textit{. }We have to show that under
assumptions (C1)-(C3), the Kuhn-Tucker conditions\textit{ }(\ref{L}) and
(\ref{CS}) are equivalent to their almost sure versions. The latter relation
is equivalent to (\ref{CS-a-s}) since the expectation of a non-negative random
variable is equal to zero if and only is the random variable is equal to zero
a.s. The former is equivalent to (\ref{l-a-s}) by virtue of (FA7).\hfill
$\square$

\section{Corollaries and applications}

\textbf{1.} \textbf{Static teams.} Consider a special case of the team model
described Sect. 2 in which all the DMs $t$ act at one moment of time $n(t)=1$
and the partial ordering on $T$ is trivial: $t\preccurlyeq s$ if and only
$t=s$. This is a \textit{static team} problem, in which all the DMs act
independently and simultaneously. Further, suppose that there are no material
constraints: $g^{t}(\omega,\cdot)=0$. Then Theorem 1 takes on the following
form (recall that we assume (A)-(F)):

\textbf{Theorem 3. }\textit{There exist stochastic Lagrange multipliers
}$q_{t}$\textit{ relaxing the informational constraints (\ref{f3-prime}), i.e.
satisfying}%
\begin{equation}
F(x)+\sum_{t\in T}Eq_{t}\left[  x_{t}-E_{t}x_{t}\right]  \leq F(\bar{x})
\label{f15}%
\end{equation}
\textit{for all feasible programs }$x=(x_{t})_{t\in T}\in X$\textit{ and each
solution }$\bar{x}$\textit{ to problem (}$\mathcal{P}$\textit{).}

Note that Theorem 3 does not involve any assumptions on the system of $\sigma
$-algebras $\mathcal{F}_{t}$. Therefore the setting at hand has the maximum
generality possible in the static case. Also, it should be noted that the
result is obtained not by the traditional method of backward induction relying
upon conditional expectations of normal integrands (see Pennanen and
Perkki\"{o} \cite{PennanenPerkkio2024}, Rockafellar and Wets
\cite{RockafellarWets1998}, and references therein). These techniques are
applicable only in the case of integral functionals and do not work in more
general situations. Although we employ, as is common, the Yosida-Hewitt
theorem (FA5) for eliminating singular functionals in $(L_{\infty})^{\ast}$ in
the construction of Lagrange multipliers, we do this in a non-standard way,
with the help of proposition (FA6) providing a possibility for a special
choice of the sequence of shrinking sets $\Gamma_{k}$ on which the singular
functional at hand is concentrated---see (\ref{f14})-(\ref{f14d}).

By applying Theorem 2 to the static team model with $g^{t}(\omega,\cdot)=0$,
we obtain

\textbf{Theorem 4.} \textit{Under assumptions (C1)-(C3), there exist
stochastic Lagrange multipliers }$q_{t}$\textit{ relaxing constraints
(\ref{f3-prime}) almost surely, i.e., satisfying for almost all }$\omega
$\textit{ and all }$a=(a_{t})\in A(\omega)$\textit{ the inequality}%

\begin{equation}
u(\omega,a)+\sum\left(  q_{t}(\omega)-\tilde{q}_{t}(\omega)\right)  a_{t}\leq
u(\omega,\bar{x}(\omega))+\sum(q_{t}(\omega)-\tilde{q}_{t}(\omega))\bar{x}%
_{t}(\omega). \label{insurance}%
\end{equation}

\textbf{2. No-regret insurance.} We briefly outline the applications of the
results obtained to insurance, referring for details to
\cite{EvstigneevKleinHaneveldandMirman1999}. By using the informational
Lagrange multipliers $q_{t}$, one can construct an insurance scheme protecting
the decision maker against uncertainty. Let us imagine that each DM $t$ pays
the \textit{premium} $\tilde{q}_{t}x_{t}$, depending on the events observed by
the DM, and gets \textit{compensation} $q_{t}x_{t}$, provided by the team
organizer after realization of all random factors. The relation $\tilde{q}%
_{t}=\mathbb{E}(q_{t}\mid\mathcal{F}_{t})$ states that the price of the
insurance policy exactly equals the expected monetary losses, i.e. we are
dealing with an \textit{actuarially fair }insurance scheme, ignoring the risk
loading (the percentage above the base premium rate). Inequality
(\ref{insurance}) shows that if every team member $t$ uses this insurance
scheme, then each component $\bar{x}_{t}$ of the optimal team program $\bar
{x}=(\bar{x}_{t})$ becomes optimal for each DM $t$ almost surely, i.e. each
agent will never regret implementing the program optimal for the whole team.
The \textit{no-regret }principle of risk bearing is systematically analyzed in
the paper \cite{EvstigneevKleinHaneveldandMirman1999} containing a whole range
of examples and applications.

\textbf{3. Shadow prices on information. }In constrained economic
optimization, the shadow price is the infinitesimal change, per infinitesimal
unit of the constraint, in the optimal value of the objective function
obtained by relaxing this constraint. In equivalent mathematical terms, this
is the corresponding Lagrange multiplier, or dual variable. The theory of
shadow prices for material constraints is well-developed, see, e.g.,
Birchenhall and Grout \cite{BirchenhallGrout1984}. Attempts to extend it to
informational ones face both technical and conceptual difficulties. The
classical Blackwell's \cite{Blackwell1951} approach, which is basically of a
discrete nature, is not suitable for dealing with infinitesimally small
increments of information. Moreover, large amounts of information understood
in the traditional Shannon's sense might have a small economic value if they
do not significantly increase the maxum value of the objective functional.

An alternative approach to measuring and evaluating information (based on the
concept of "flexibility" of a decision \ function) is proposed in
\cite{EvstigneevKleinHaneveldandMirman1999}. It fits well the current context:
one can show that the informational Lagrange multipliers $q_{t}$ in
(\ref{insurance}) represent the shadow prices on information defined in the
usual, infinitesimal terms. A comprehensive analysis of these prices is
conducted in \cite{EvstigneevKleinHaneveldandMirman1999}. However, the model
in \cite{EvstigneevKleinHaneveldandMirman1999} deals with a single DM. It
would be of interest to develop an analogous theory for teams with multiple DMs.

\section{Appendix: General facts from functional analysis}

\textbf{1.} Let $D_{1}$ and $D_{2}$ be Banach spaces, $X$ a convex subset in
$D_{1}$, $F$ a concave lower-bounded functional defined on $X$ and
$h:D_{1}\rightarrow D_{2}$ a continuous linear mapping. Assume that $F(x)$
attains its maximum on the set $\{x:x\in X,h(x)=0\}$ at some point $\bar{x}$.

\begin{itemize}
\item[\textbf{(FA1)}] If $\{x:h(x)=0\}$ contains a point $\mathring{x}$
belonging to the interior Int $X$ of the set $X$ and $h$ maps $D_{1}$ onto a
closed subspace $D_{3}\subset D_{2}$, then there exists a functional
$\varphi\in D_{2}^{\ast}$ ($D_{2}^{\ast}$ being the dual of $D_{2}$) such
that
\[
F(x)+\varphi(h(x))\leq F(\bar{x}),\ \ x\in X.
\]
(An infinite-dimensional version of the Kuhn-Tucker theorem for equality
constraints.)\hfill$\square$
\end{itemize}

\textit{Proof.} Let $U$ be an open subset of $X$ such that $\mathring{x}\in
U$. Since $h$ maps $D_{1}$ onto $D_{3}$, then by the Open Mapping Theorem
(see, e.g. \cite{DunfordSchwartz1964}, II.2.1), $h(U)$ is open in $D_{3}$.
Thus the set
\[
A=\{(r,y):r\leq F(x),y=h(x),x\in X\}
\]
contains the open set
\[
A^{\prime}=\{(r,y):r<c,y\in h(U)\}\subset\mathbb{R}^{1}\times D_{3},
\]
where $c$ is a constant such that $F(x)\geq c$ for $x\in X$. Consider the set
\[
B=\{(r,y):r>F(\bar{x}),\ y=0\}.
\]
Since $\operatorname{Int}A\neq\varnothing$ and $A\cap B=\varnothing$, the sets
$A$ and $B$ can be separated, i.e., there is a non-zero linear functional
$\lambda=(\lambda_{1},\lambda_{2})$ on $\mathbb{R}^{1}\times D_{3}$ such that
$\lambda(a)\leq\lambda(b)$ for $a\in A$, $b\in B$ (see, e.g.,
\cite{DunfordSchwartz1964}, V.3.8). Now we have
\[
\lambda_{1}F(x)+\lambda_{2}h(x)\leq\lambda_{1}F(\bar{x})\quad(x\in X).
\]
Hence, if $\lambda_{1}=0$, then $\lambda_{2}(y)\leq0$ for $y\in h(U)$, where
$h(U)$ is open and contains $0$. Thus $\lambda_{2}=0$ in contrary to the fact
that $(\lambda_{1},\lambda_{2})\neq0$. Since $\lambda_{1}\neq0$, we can define
a functional $\varphi^{\prime}:=\lambda_{2}/\lambda_{1}$. In view of the
Hahn-Banach theorem (\cite{DunfordSchwartz1964}, II.3.11), there exists a
linear functional $\varphi\in D_{2}^{\ast}$ such that $\varphi=\varphi
^{\prime}$ on $D_{3}$. The functional $\varphi$ has the desired
property.\hfill$\square$

\textbf{2.} Suppose now that the space $D_{2}$ is partially ordered. Namely, a
convex cone $K_{2}\subset D_{2}$ is fixed and, by the definition, $y\geq x$ if
$y-x\in K_{2}$. Let $g:X\rightarrow D_{2}$ be a concave mapping. (The
concavity of a mapping is defined exactly in the same way as the concavity of
a real-valued function.)

\begin{itemize}
\item[\textbf{(FA2)}] Assume that the maximum of $F$ on $\{x:x\in
X,g(x)\geq0\}$ is attained at some $\tilde{x}$ and there is $\hat{x}$ with
$g(\hat{x})\in\operatorname{Int}K_{2}$. Then one can find a functional $\pi\in
D_{2}^{\ast}$ with properties: $\pi\geq0$ (i.e., $\pi(x)\geq0$ for $x\geq0$)
and $F(x)+\pi(g(x))\leq F(\tilde{x})$, $x\in X$.
\end{itemize}

This is an infinite-dimensional version of the Kuhn-Tucker theorem for concave
constraints of inequality type. A proof can be found, e.g., in Hurwicz
\cite{Hurwicz1958}.

\textbf{3.} We shall now formulate some results concerning the Banach spaces
$L_{1}$ and $L_{\infty}$. The space $L_{1}=L_{1}(\mathcal{F},V)$ consists of
integrable $\mathcal{F}$-measurable vector-functions with values in $V$ and
$L_{\infty}=L_{\infty}(\mathcal{F},V)$ is the space of essentially bounded
$\mathcal{F}$-measurable vector-functions with values in $V$. The norms of
$L_{1}$ and $L_{\infty}$ will be denoted, respectively, by $\Vert\cdot
\Vert_{1}$ and $\Vert\cdot\Vert_{\infty}$.

\begin{itemize}
\item[\textbf{(FA3)}] (Koml\'{o}s \cite{Komlos1967}). Let $x_{1},x_{2},\ldots$
be a sequence of elements of $L_{1}$. If $\sup\Vert x_{i}\Vert_{1}<\infty$,
then there exists a sequence of natural numbers $m_{1},m_{2},...$ and $x\in
L_{1}$ such that
\begin{equation}
\frac{x_{m_{1}}+x_{m_{2}}+\cdots+x_{m_{k}}}{k}\rightarrow x\ \text{a.s.},
\label{f6}%
\end{equation}
and (\ref{f6}) holds for any subsequence of the sequence $m_{1},m_{2},...$ $.$

\item[\textbf{(FA4)}] Let $Q$ be a convex bounded subset of $L_{1}$ and $F$ a
concave functional on $Q$. If $Q$ is closed and $F$ is upper semi-continuous
with respect to the convergence almost everywhere, then $F$ attains its
maximum on $Q$.
\end{itemize}

\textit{Proof.} Consider a sequence $x_{m}\in Q$ such that
\[
F(x_{m})\rightarrow\bar{F}=\sup\{F(x):x\in Q\}.
\]
By (FA3), for some $(m_{k})$ and $x\in Q$, we have
\[
F(x)\geq\overline{\lim}F\left(  \frac{x_{m_{1}}+\cdots+x_{m_{k}}}{k}\right)
\geq\lim\frac{1}{k}(F(x_{m_{1}})+\cdots+F(x_{m_{k}}))=\bar{F}.
\]
Thus, $F(x)=\bar{F}$ i.e., $x$ maximizes $F$ on $Q$. \hfill$\square$

\textbf{4.} A functional $\pi\in L_{\infty}^{\ast}(\mathcal{F},R^{n})$ is
called \textit{absolutely continuous} if it can be represented in the form
$\pi(x)=Ex_{0}x$ with $x_{0}\in L_{1}$. A functional $\pi\in L_{\infty}^{\ast
}$ is called \textit{singular} if there is a sequence of sets $\Gamma
_{1},\Gamma_{2},\ldots$ such that
\begin{equation}
\pi(x)=\pi(\chi_{\Gamma_{k}}x)\quad(\forall k)\quad\text{and}\quad
P(\Gamma_{k})\rightarrow0, \label{f7}%
\end{equation}
where $\chi_{\Gamma}=1$ on $\Gamma$ and $\chi_{\Gamma}=0$ on the complement
$\bar{\Gamma}$ of $\Gamma$.

\begin{itemize}
\item[\textbf{(FA5)}] (Yosida and Hewitt \cite{YosidaHewitt1952}). Every
functional $\pi\in L_{\infty}^{\ast}(\mathcal{F},R^{n})$ can be uniquely
decomposed into the sum $\pi=\pi^{a}+\pi^{s}$ where $\pi^{a}$ is absolutely
continuous and $\pi^{s}$ is singular. If $\pi\geq0$, then $\pi^{a}\geq0$ and
$\pi^{s}\geq0$.

\item[\textbf{(FA6)}] Let $\pi$ be a singular functional on $L_{\infty
}=L_{\infty}(\mathcal{F},V)$ and $\mathcal{G}$ a sub-$\sigma$-algebra of
$\mathcal{F}$. There exist sets $\Gamma_{k}$ possessing the properties
(\ref{f7}) and, in addition, the following property:
\begin{equation}
\Vert\chi_{\bar{\Gamma}_{k}}P(\Gamma_{k}\mid\mathcal{G})\Vert_{\infty
}\rightarrow0, \label{f8}%
\end{equation}
where $\bar{\Gamma}_{k}=\Omega\backslash\Gamma_{k}$.
\end{itemize}

\textit{Proof.} Let $\Delta_{k}$ be any sequence of sets with the property
(\ref{f7}). We may assume without loss of generality that $P(\Delta_{k})\leq
k^{-2}$. Consider the sets
\[
\Delta_{k}^{\prime}=\{\omega:P(\Delta_{k}\mid\mathcal{G})\geq\frac{1}%
{k}\}\quad\text{and}\quad\Gamma_{k}=\Delta_{k}\cup\Delta_{k}^{\prime}.
\]
Since $\Gamma_{k}\supseteq\Delta_{k}$, we have $\pi(\chi_{\Gamma_{k}}%
x)=\pi(x)$. Further, $P(\Gamma_{k})\rightarrow0$ because
\[
P(\Gamma_{k})\leq P(\Delta_{k})+P(\Delta_{k}^{\prime-2}+k^{-1}.
\]
\smallskip Finally, (\ref{f7}) follows from the relations:
\[
\Vert\chi_{\bar{\Gamma}_{k}}E(\chi_{\Gamma_{k}}\mid\mathcal{G})\Vert_{\infty
}\leq\Vert\chi_{\bar{\Delta}_{k}}\chi_{\bar{\Delta}_{k}^{\prime}}E\left(
(\chi_{\Delta_{k}}+\chi_{\bar{\Delta}_{k}^{\prime}})\mid\mathcal{G})\right)
\Vert_{\infty}%
\]%
\[
=\Vert\chi_{\bar{\Delta}_{k}}\chi_{\bar{\Delta}_{k}^{\prime}}E\left(
(\chi_{\Delta_{k}}\mid\mathcal{G})\right)  \Vert_{\infty}\leq\frac{1}{k}.
\]

\hfill$\square$

(The equality in the above formula is valid since $E(\chi_{\Delta_{k}^{\prime
}}\mid\mathcal{G})=\chi_{\Delta_{k}^{\prime}}$.)

\textbf{5.} Let $(\Omega,\mathcal{F},P)$ be a probability space, $A(\omega)$ a
compact subset of a finite-dimensional space and $f(\omega,a)$ ($\omega
\in\Omega,a\in A(\omega)$) a real function.

\begin{itemize}
\item[\textbf{(FA7)}] If $A(\omega)$ and $f(\omega,a)$ have the properties
(C1)--(C3), then there exists an $\mathcal{F}$-measurable function $\tilde
{x}(\omega)$ such that
\begin{equation}
f(\omega,\tilde{x}(\omega))=\max_{a\in A(\omega)}f(\omega,a)\ \text{a.s.},
\label{f-max}%
\end{equation}
and relation (\ref{f-max}) holds if and only if
\[
Ef(\omega,x(\omega))\leq Ef(\omega,\tilde{x}(\omega))
\]
for all $x(\omega)$ satisfying $x(\omega)\in A(\omega)$ (a.s.).

\item[\textbf{(FA8)}] Let $\mathcal{G}$ be a sub-$\sigma$-algebra of
$\mathcal{F}$ and $x(\omega)$ an $\mathcal{F}$-measurable integrable function.
Assume that $A(\omega)$ depends $\mathcal{G}$-measurably on $\omega$ and is
convex for each $\omega$. If $x(\omega)\in A(\omega)$ a.s., then
$E(x(\omega)\mid\mathcal{G})\in A(\omega)$ a.s.
\end{itemize}

For a proof of (FA8) see, Arkin and Evstigneev \cite{ArkinEvstigneev1987},
Appendix II, Lemma 1.

\end{document}